\documentclass[11pt,leqno]{article}
\usepackage{exscale,relsize}
\usepackage{amsmath}
\usepackage{amsfonts}
\usepackage{amssymb}
\usepackage{calc}
\usepackage{theorem}
\usepackage{pifont}      
\newcommand{\scal}[2]{\langle{{#1},{#2}}\rangle}

\newcommand{\RR}{\ensuremath{\mathbb R}}

\newcommand{\RX}{\ensuremath{\,\left]-\infty,+\infty\right]}}
\newcommand{\RXX}{\ensuremath{\,\left[-\infty,+\infty\right]}}

\newcommand{\menge}[2]{\big\{{#1} \mid {#2}\big\}}

\newcommand{\aff}{\operatorname{aff}}

\newcommand{\cone}{\ensuremath{\operatorname{cone}}}
\newcommand{\dom}{\ensuremath{\operatorname{dom}}}

\newcommand{\gra}{\ensuremath{\operatorname{gra}}}

\newcommand{\conv}{\ensuremath{\operatorname{conv}}}

\renewcommand{\phi}{\ensuremath{\varphi}}

\newcommand{\To}{\ensuremath{\rightrightarrows}}

\newtheorem{theorem}{Theorem}[section]

\newtheorem{fact}[theorem]{Fact}
\newtheorem{corollary}[theorem]{Corollary}

\newtheorem{definition}[theorem]{Definition}

\theoremstyle{plain}{\theorembodyfont{\rmfamily}
}
\theoremstyle{plain}{\theorembodyfont{\rmfamily}
}
\theoremstyle{plain}{\theorembodyfont{\rmfamily}
}
\theoremstyle{plain}{\theorembodyfont{\rmfamily}
}
\theoremstyle{plain}{\theorembodyfont{\rmfamily}
\newtheorem{remark}[theorem]{Remark}}
\theoremstyle{plain}{\theorembodyfont{\rmfamily}
}

\begin{document}

\title{{\sffamily An affirmative answer to a problem posed by Z\u{a}linescu}}

\author{
Liangjin\
Yao\thanks{Mathematics, Irving K.\ Barber School, UBC Okanagan,
Kelowna, British Columbia V1V 1V7, Canada.
E-mail:  \texttt{ljinyao@interchange.ubc.ca}.}. }

\date{November 9, 2009}

\maketitle


\begin{abstract} \noindent
Recently, in \cite{Zali} Z{\v{a}}linescu posed a question about the characterization
 of the intrinsic core of the Minkowski sum of two graphs
associated with two maximal monotone operators.
In this note we give an affirmative answer.
\end{abstract}



\section{Introduction}
We suppose throughout this note that
$X$ is  a real reflexive Banach space with norm $\|\cdot\|$
and dual product $\scal{\cdot}{\cdot}$. We now introduce some notation.
Let $A\colon X\To X^*$
be a \emph{set-valued operator} or \emph{multifunction} whose graph is defined by
\begin{align*}\gra A:= \menge{(x,x^*)\in X\times X^*}{x^*\in Ax}.\end{align*}
The \emph{domain} of $A$ is $\dom A := \menge{x\in X}{Ax\neq\varnothing}$.
 Recall that $A$ is  \emph{monotone} if for all $(x,x^*), (y,y^*)\in\gra
A$ we have
\begin{equation*}
\scal{x-y}{x^*-y^*}\geq 0,
\end{equation*}and $A$ is \emph{maximal monotone} if $A$ is monotone and
 $A$ has no proper monotone extension (in the sense of graph inclusions).

The
\emph{Fitzpatrick function} of $A$ (see \cite{Fitz88}) is given by
\begin{equation}
F_A\colon  (x,x^*)\mapsto \sup_{(a,a^*)\in
\gra A}\big(\langle x, a^*\rangle+ \langle a,x^*\rangle-\langle a,a^*\rangle\big).
\end{equation}
For a function  $f\colon X\to \RX$, the \emph{domain} is
$\dom f:= \{x\in X \mid f(x)<+\infty\}$ and
$f^*\colon X^*\to\RXX\colon x^*\mapsto\sup_{x\in X}(\scal{x}{x^*}-f(x))$ is
the \emph{Fenchel conjugate} of $f$.

Given $F:X\times X^*\rightarrow\RX$, we say $F$ is a \emph{representative} of a maximal monotone
operator $A$ if $F$ is lower semicontinuous and convex with
 $F\geq \langle\cdot,\cdot\rangle$, $F^*\geq \langle\cdot,\cdot\rangle$
and
\begin{align*}
\gra A=\{(x,x^*)\mid F(x,x^*)=\langle x,x^*\rangle\}.\end{align*}

Following \cite{Penot2}, it will be convenient to set
$ F^\intercal\colon X^*\times X:\rightarrow\RX\colon (x^*,x)\mapsto F(x,x^*)$,
where
$F\colon X\times X^*\to\RX$, and similarly for a function defined on $X^*\times X$.

We define $\widehat{F}$ (see \cite{Zali}) by
\begin{align*}
\widehat{F}(x,x^*):=F(x,-x^*).\end{align*}
Let $a=(x,x^*), b=(y,y^*)\in X\times X^*$, we also set (see \cite{Si2}) by
\begin{align*}
\lfloor a,b\rfloor=\langle x,y^*\rangle+\langle y,x^*\rangle.\end{align*}
Given a subset $D$ of $X$,
$\overline{D}$ is the \emph{closure},
$\conv{D}$ is the \emph{convex hull}, and
$\aff{D}$ is the \emph{affine hull}. The \emph{conic hull} of $D$ is denoted by
$\cone D:=\{\lambda x\mid \lambda\geq0, x\in D\}$.
The \emph{indicator function} $\iota_D:X\rightarrow\RX$ of $D$ is defined by
\begin{equation*}x\mapsto\begin{cases}0,\;&\text{if}\;x\in D;\\
+\infty,\;&\text{otherwise}.\end{cases}\end{equation*}
The \emph{intrinsic core} or \emph{relative algebraic interior} of $D$, written as $^{i}D$ in \cite{Zalinescu}, is
\begin{align*}^{i}D:=\{a\in D\mid \forall x\in \aff(D-D),
\exists\delta>0, \forall\lambda\in\left[0,\delta\right]:
a+\lambda x\in D\}.
\end{align*}
We define $^{ic}D$ by
\begin{equation*}
^{ic}D:=\begin{cases}^{i}D,\,&\text{if $\aff D$ is closed};\\
\varnothing,\,&\text{otherwise}.
\end{cases}
\end{equation*}

 Z{\v{a}}linescu posed the following problem in \cite{Zali}:
Let $A, B:X\rightrightarrows X^*$ be  maximal monotone. Is the implication
\begin{align*}
^{ic}\left[\conv{(\gra A-\gra (-B))}\right]\neq\varnothing\quad\Rightarrow\quad
^{ic}\left[\dom F_A-\dom \widehat{F_B}\right]\neq\varnothing\end{align*}
true?
Theorem~\ref{Zopen} provides an affirmative answer to this question. It further shows that
these two sets actually are equal.

\section{Main result}

\begin{definition}[Fitzpatrick family]
Let  $A\colon X \To X^*$ be a maximal monotone operator.
The associated \emph{Fitzpatrick family}
$\mathcal{F}_A$ consists of all functions $F\colon
X\times X^*\to\RX$ that are lower semicontinuous and convex,
and that satisfy
$F\geq \scal{\cdot}{\cdot} $, and $F=\scal{\cdot}{\cdot}$ on $\gra A$.
\end{definition}

\begin{fact}[Fitzpatrick]\label{GF:2}
\emph{(See \cite[Theorem~3.10]{Fitz88}.)}
Let  $A\colon X \To X^*$ be a maximal monotone operator.
Then for every $(x,x^*)\in X\times X^*$,
\begin{equation}
F_A(x,x^*) = \min\menge{F(x,x^*)}{F\in \mathcal{F}_A}
\quad\text{and}\quad
F_A^{*\intercal}(x,x^*) = \max\menge{F(x,x^*)}{F\in \mathcal{F}_A}.
\end{equation}
\end{fact}

\begin{fact}[Simons]\emph{(See \cite[Lemma~20.4(b)]{Si2}.)}\label{Sf:1}
Let $A:X\rightrightarrows X^*$ be maximal monotone and $a:=(x,x^*)\in X\times X^*$ with
$\langle x,x^*\rangle=0$.
 Suppose that there exists $u\in\RR$ such that
\begin{align*}
\lfloor \gra A,a\rfloor=\{u\}.\end{align*}
Then
\begin{align*}
\lfloor \dom F_A,a\rfloor=\{u\}.\end{align*}

\end{fact}

\begin{theorem}\label{zad}
Let $A, B:X\rightrightarrows X^*$ be  maximal monotone. Then
\begin{align}
\overline{\aff\left[\gra A-\gra(-B)\right]}&=\overline{\aff\left[\dom F_A-\dom \widehat{F_B}\right]}\label{Sf:2}
.\end{align}

\end{theorem}

\begin{proof} We do and can suppose $(0,0)\in \gra A$ and $(0,0)\in\gra B$.
We first show
\begin{align}\left[\dom F_A-\dom \widehat{F_B}\right]\subseteq\overline{\aff\left[\gra A-\gra(-B)\right]}.\end{align}
Suppose to the contrary that there exists $c\in X\times X^*$ such that
$c\in \left[\dom F_A-\dom \widehat{F_B}\right]$ but $c\notin\overline{\aff\left[\gra A-\gra (-B)\right]}$.
By the Separation Theorem, there exist $a:=(x,x^*)\in X\times X^*$ and $\delta\in\RR$ such that
\begin{align}
\lfloor a, c\rfloor>\delta>\sup\big\{\lfloor a, e\rfloor\mid e\in\overline{\aff\left[\gra A-\gra (-B)\right]}\big\}.\label{TF:1}\end{align}
Since $(0,0)\in\gra A$, $(0,0)\in\gra B$ and
$\overline{\aff\left[\gra A-\gra(-B)\right]}$ is a closed subspace, we have $\delta>0$ and
$\lfloor a, b-d\rfloor=0, \forall b\in\gra A,\, \forall d\in\gra(-B)$.
Thus,
\begin{align}
\lfloor a, \gra A\rfloor=\{0\}=\lfloor a,\gra (-B)\rfloor=\lfloor (-x,x^*),\, \gra B\rfloor.\label{TF:2}\end{align}
By $(0,0)\in\gra A$ and $(0,0)\in\gra B$ again,
\begin{align} F_A(a)=F_A(x,x^*)=0, \quad F_B(-x, x^*)=0.\label{TF:3}\end{align}
Since $F_A(x,x^*)\geq\langle x,x^*\rangle $ and $F_B(-x, x^*)\geq\langle-x,x^*\rangle$,
by \eqref{TF:3},
$\langle x, x^*\rangle=0$.
Thus by \eqref{TF:2} and Fact~\ref{Sf:1},
\begin{align*}
\lfloor a, \dom F_A\rfloor=\{0\}=\lfloor (-x,x^*),\,\dom F_B\rfloor
=\lfloor a, \dom \widehat{F_B}\rfloor.\end{align*}
Thus,
$\lfloor a, \dom F_A-\dom \widehat{F_B}\rfloor=\{0\}$, which
contradicts \eqref{TF:1}.
Hence \begin{align}\left[\dom F_A-\dom \widehat{F_B}\right]\subseteq\overline{\aff\left[\gra A-\gra(-B)\right]}.\end{align}
And thus $\overline{\aff\left[\dom F_A-\dom \widehat{F_B}\right]}\subseteq\overline{\aff\left[\gra A-\gra(-B)\right]}$.
Hence
\begin{align*}
\overline{\aff\left[\gra A-\gra(-B)\right]}&=\overline{\aff\left[\dom F_A-\dom \widehat{F_B}\right]}.
\end{align*}
\end{proof}

\begin{fact}[Z{\v{a}}linescu]\emph{(See \cite[Lemma~2 and Theorem~3]{Zali}.)}\label{SFF:1}
Let $A, B:X\rightrightarrows X^*$ be  maximal monotone, and let $F_1, F_2$ be  representatives
 of $A, B$, respectively.   Then
 \begin{align*}{^{ic}\left[\gra A-\gra (-B)\right]}={^{ic}\left[\conv(\gra A-\gra (-B))\right]}\end{align*}
 and
\begin{align*}
{^{ic}\left[\dom F_1 -\dom \widehat{F_2}\right]}\subseteq\left[\gra A-\gra (-B)\right]
\subseteq \conv\left[\gra A-\gra (-B)\right]\subseteq\left[\dom F_1 -\dom \widehat{F_2}\right]
.\end{align*}
If ${^{ic}\left[\dom F_1 -\dom \widehat{F_2}\right]}\neq\varnothing$, then  \begin{align}
^{ic}\left[\dom F_1 -\dom \widehat{F_2}\right]={^{ic}\left[\gra A-\gra (-B)\right]}
= {^{ic}\left[\conv(\gra A-\gra (-B))\right]}
.\label{SF:19}\end{align}
\end{fact}

\begin{remark}
If $X$ is finite-dimensional, the intrinsic core of a convex set $D\subseteq X$ is the same as the
relative interior of $D$ in the sense of Rockafellar\cite{Rocky}. Then ${^{ic}\left[\dom F_1 -\dom \widehat{F_2}\right]}
= {^{i}\left[\dom F_1 -\dom \widehat{F_2}\right]}\neq\varnothing$ by \cite[Theorem~6.2]{Rocky}.
Thus, \eqref{SF:19} always holds.
\end{remark}

Our main result comes the following which  provides an affirmative answer to the question
 posed by Z\u{a}linescu.

\begin{theorem}\label{Zopen}
Let $A, B:X\rightrightarrows X^*$ be  maximal monotone such that
${^{ic}\left[\conv{(\gra A-\gra (-B))}\right]}\neq\varnothing$. Then
\begin{align*}{^{ic}\left[\gra A-\gra (-B)\right]}=
{^{ic}\left[\conv{(\gra A-\gra (-B))}\right]}={^{ic}\left[\dom F_A-\dom \widehat{F_B}\right]}.
\end{align*}
Moreover, if $F_1, F_2$ are  representatives
 of $A, B$, respectively, then
\begin{align}{^{ic}\left[\dom F_1 -\dom \widehat{F_2}\right]}={^{ic}\left[\conv{(\gra A-\gra (-B))}\right]}
={^{ic}\left[\dom F_A-\dom \widehat{F_B}\right]}.\label{SF:4}\end{align}
\end{theorem}
\begin{proof}Let $a\in{^{ic}\left[\conv{(\gra A-\gra (-B))}\right]}$. Then we have $a\in\left[\dom F_A-\dom \widehat{F_B}\right]$
and $\cone \left[\conv{(\gra A-\gra (-B))}-a\right]$ is a closed subspace.
By Theorem~\ref{zad},
\begin{align*}
&\cone \left[\conv{(\gra A-\gra (-B))}-a\right]\subseteq\cone \left[\dom F_A-\dom \widehat{F_B}-a\right]\\
&\subseteq\aff \left[\dom F_A-\dom \widehat{F_B}-a\right]=\aff\left[\dom F_A-\dom \widehat{F_B}\right]-\{a\}\\
&\subseteq\overline{\aff \left[\gra A-\gra (-B)\right]}-\{a\}
\subseteq\overline{\aff \left[\gra A-\gra (-B)-a\right]}\\&=\overline{\aff \left[\conv{(\gra A-\gra (-B))}-a\right]}
\subseteq \cone \left[\conv{(\gra A-\gra (-B))}-a\right].
\end{align*}
Hence $\cone \left[\dom F_A-\dom \widehat{F_B}-a\right]=\cone \left[\conv{(\gra A-\gra (-B))}-a\right]$
is a closed subspace.
Thus
$a\in{^{ic}\left[\dom F_A-\dom \widehat{F_B}\right]}$. By Fact~\ref{SFF:1},
\begin{align*}^{ic}\left[\gra A-\gra (-B)\right]=
{^{ic}\left[\conv{(\gra A-\gra (-B))}\right]}={^{ic}\left[\dom F_A-\dom \widehat{F_B}\right]}.
\end{align*}
And by Fact~\ref{GF:2},
\begin{align*}
\conv{\left[\gra A-\gra (-B))\right]}
\subseteq
\left[\dom F_1 -\dom \widehat{F_2}\right]\subseteq\left[\dom F_A-\dom \widehat{F_B}\right]
.\end{align*}
Similar to the proof above,  see that \eqref{SF:4} holds.
\end{proof}

\begin{theorem}\label{Zopen:1}
Let $A, B:X\rightrightarrows X^*$ be  maximal monotone, and $F_1, F_2$ be  representatives
 of $A, B$, respectively.  Then
\begin{align}{^{ic}\left[\dom F_1 -\dom \widehat{F_2}\right]}={^{ic}\left[\conv{(\gra A-\gra (-B))}\right]}
={^{ic}\left[\gra A-\gra (-B)\right]}.\label{SF:5}\end{align}
\end{theorem}
\begin{proof}
We consider two cases.

Case 1: ${^{ic}\left[\conv{(\gra A-\gra (-B))}\right]}=\varnothing$. Assume
 that  ${^{ic}\left[\dom F_1 -\dom \widehat{F_2}\right]}\neq\varnothing$.
Then by
Fact~\ref{SFF:1}, ${^{ic}\left[\dom F_1 -\dom \widehat{F_2}\right]}={^{ic}\left[\conv{(\gra A-\gra (-B))}\right]}=\varnothing$.
This a contradiction.

Case 2: ${^{ic}\left[\conv{(\gra A-\gra (-B))}\right]}\neq\varnothing$. Apply Theorem~\ref{Zopen}.

Combining the above results, we see that \eqref{SF:5} holds.

\end{proof}

\begin{corollary}
Let $A, B:X\rightrightarrows X^*$ be  maximal monotone, and $F_1, F_2$ be representatives
 of $A, B$, respectively. Assume ${^{ic}\left[\conv{(\gra A-\gra (-B))}\right]}\neq\varnothing$. Then
\begin{align}\overline{\left[\dom F_1 -\dom \widehat{F_2}\right]}=\overline{\conv{\left[(\gra A-\gra (-B)\right]}}
=\overline{\left[\gra A-\gra (-B)\right]}.\label{SF:12}\end{align}
In particular,
\begin{align*}\overline{\left[\dom F_A -\dom \widehat{F_B}\right]}=\overline{\conv{\left[(\gra A-\gra (-B)\right]}}
=\overline{\left[\gra A-\gra (-B)\right]}.\end{align*}
\end{corollary}

\begin{proof}
Given a convex set $D\subseteq X$, assume that ${^{ic}D}\neq\varnothing$, then $\overline{^{ic}D}=\overline{D}$.
By Theorem~\ref{Zopen:1},
\begin{align*}
\overline{\left[\dom F_1 -\dom \widehat{F_2}\right]}=\overline{\conv{\left[(\gra A-\gra (-B)\right]}}
=\overline{^{ic}\left[\gra A-\gra (-B)\right]}\subseteq\overline{\left[\gra A-\gra (-B)\right]}\end{align*}
Hence \eqref{SF:12} holds.
\end{proof}

\section*{Acknowledgment} The author thanks
Heinz Bauschke and  Xianfu Wang for valuable discussions.



\end{document}